\documentclass[reqno, 12pt]{amsart}
\usepackage{graphicx}
\usepackage{psfrag}
\usepackage{epsfig,color}
\input{psfig.sty}
\usepackage{amsmath}
\usepackage{graphicx}
\newcommand{\be}{\begin{eqnarray}}
\newcommand{\ee}{\end{eqnarray}}
\newcommand{\ben}{\begin{eqnarray*}}
\newcommand{\een}{\end{eqnarray*}}
\newtheorem{teo}{Theorem}[section]

\newtheorem{remark}[teo]{Remark}

\newtheorem{proposition}[teo]{Proposition}

\newtheorem{theorem}[teo]{Theorem}
\newtheorem{lemma}[teo]{Lemma}
\usepackage{amsbsy}
\newcommand\RR{{\mathbb R}}

\numberwithin{equation}{section}

\begin{document}
\title[On Regions of Existence and Nonexistence]
{On Regions of Existence and Nonexistence of solutions for a System of $p$-$q$-Laplacians}
\thanks{This research was supported by
        FONDECYT 7040105 for the first author, FONDECYT 1030593 for the second author,
        Fondap Matem\'aticas Aplicadas and FONDECYT 3040059 for the third author,
        Fondap Matem\'aticas Aplicadas and Milenio
        Grant P01-34 for the fourth author.}
\author[Ph. Cl\'ement]{\qquad\quad  Philippe Cl\'ement}
\address{ Department of Technical Mathematics and Informatics,
University of Technology, Delft,
The Netherlands}
\email{\tt ph.p.j.e.clement@ewi.tudelft.nl}

\author[M. Garc\'\i a-Huidobro]{\qquad  Marta Garc\'{\i}a-Huidobro}
\address{Departamento de Matem\'atica,
        Universidad Cat\'olica de Chile,
        Casilla 306, Correo 22,
        Santiago, Chile.}
\email{\tt mgarcia@mat.puc.cl}

\author[I. Guerra]{\qquad  Ignacio Guerra\quad}
\address{ Centro de Modelamiento Matemático,
 Universidad  de Chile,
Casilla  170 Correo 3, Santiago, Chile}
\email{\tt iguerra@dim.uchile.cl}

\author[R. Man\'asevich]{\quad\quad  Ra\'ul Man\'asevich}
\address{Departamento de Ingenier\'\i a Matem\'atica, FCFM, Universidad de Chile,
Casilla 170, Correo 3,
Santiago, Chile}
\email{manasevi@dim.uchile.cl}

\subjclass{Primary 35 J 70, Secondary 35 J 60}
\keywords{$m$-Laplacian, energy identities.  }


\begin{abstract}
We give a new region of existence of  solutions to the superhomogeneous
Dirichlet problem
$$ \quad
\begin{array}{l}
-\Delta_{p} u= v^\delta\quad  v>0\quad \mbox{in}\quad B,\cr
-\Delta_{q} v = u^{\mu }\quad  u>0\quad \mbox{in}\quad B, \cr u=v=0
\quad \mbox{on}\quad \partial B,
\end{array}\leqno{(S_R)}
$$
where $B$ is the ball of radius $R>0$ centered at the origin in
$\RR^N.$ Here $\delta, \mu >0$ and
$
\Delta_{m} u={\rm div}(|\nabla u|^{m-2}\nabla u)
$
is the $m-$Laplacian operator for $m>1$.

\end{abstract}

\maketitle

\section{Introduction and main results}
Consider the quasilinear elliptic system
$$ \quad
\begin{array}{l}
-\Delta_{p} u= v^\delta\quad  v>0\quad \mbox{in}\quad B,\cr
-\Delta_{q} v = u^{\mu }\quad  u>0\quad \mbox{in}\quad B, \cr
u=v=0 \quad \mbox{on}\quad \partial B,
\end{array}\leqno{(S_R)}
$$
where $B$ is the ball of radius $R>0$ centered at the origin in $\RR^N.$ Here
$\delta, \mu >0$ and
$$
\Delta_{m} u={\rm div}(|\nabla u|^{m-2}\nabla u)
$$
is the $m-$Laplacian operator for $m>1$.

In view of the invariance of problem $(S_R)$ under rotations, it
is natural to look for radially symmetric solutions. If we still
denote by $u, v$ the solutions as functions of $r=|x|$, we obtain
the system of ode's
\be\label{1.1}
\begin{gathered}
-(r^{N-1}|u'(r)|^{p-2}u'(r))'=r^{N-1}|v(r)|^\delta,\quad 0<r<R\\
-(r^{N-1}|v'(r)|^{q-2}v'(r))'=r^{N-1}|u(r)|^\mu ,\quad 0<r<R,
\end{gathered}
\ee with appropriate boundary conditions. We are primarily
interested in the existence of (regular) solutions of \eqref{1.1},
i.e.,  $(u,v)\in (C^1[0,R]\cap C^2(0,R])^2$ satisfying \eqref{1.1}
and $u'(0)=v'(0)=0$, $u(R)=v(R)=0$.

Clearly, either both $u$ and $v$ are identically $0$, or both $u$ and $v$ are strictly positive
and decreasing on $[0,R)$.

Observe that system $(S_R)$ is homogeneous in the sense that if
$(u,v)$ is a solution, then $(\lambda u,\nu v)$ is also a
solution provided that $\lambda,\ \nu>0$ and
$\lambda^{1-p}=\nu^\delta$ and $\nu^{1-q}=\lambda^\mu $. So it is
natural to call the system {\em superhomogeneous} when
$$
 d:=\delta \mu -(p-1)(q-1)>0\quad \delta>0,\:\mu >0.\leqno(H_1)
$$
In case that $p=q=2$, condition $(H_1)$ is usually called {\em
superlinear condition} and it is equivalent to the condition
$$\frac{1}{\delta+1}+\frac{1}{\mu +1}<1.$$

It has been shown in \cite{CFM92}, \cite{Mi93}, \cite{PevdV92}, and  \cite{V} that under $(H_1)$, when $N>2$, a necessary and
sufficient condition for the existence of   radial solutions to $(S_R)$ is
$$\frac{1}{\delta+1}+\frac{1}{\mu +1}>\frac{N-2}{N}.$$

In case that $m=p=q\not=2$ and $\delta=\mu $ (see Remark
\ref{poneeq} in the appendix) we have that if $(u,v)$ is a
solution, then $u=v$ and hence the system reduces to an equation.
It follows then from results of \cite{O} that a   solution exists
in that case (for $m<N$) if and only if \be\label{eqcond}
\frac{1}{\delta+1}>\frac{N-m}{Nm}. \ee Apart from these cases no
necessary and sufficient condition for the existence of solutions
is known. Sufficient conditions have been obtained in \cite{CMM93}
where a-priori estimates are established by means of a blow up
method in the sense of Gidas and Spruck, see \cite{G-S}, and a
degree argument. In \cite{ACM}, the problem has been studied in a
bounded convex domain
with $C^2$ boundary.

The main goal of this paper is to exhibit a new region of
existence of   solutions to $(S_R)$.
 This is done in Theorem \ref{mainth2}.

To our knowledge, when $p\not=q$ or $p=q\not=2$ and $\delta\not=\mu $, there are no nonexistence results
(of Pohozaev type) in
the literature. In Theorem \ref{mainth} we provide such a region of nonexistence.

An important ingredient in the proof of our main result Theorem
\ref{mainth2} is the observation that under condition $(H_1)$, the
absence of positive \lq\lq ground states\rq\rq\ implies existence
of solutions for $(S_R)$. The result is contained implicitly in
\cite{CMM93}, \cite{CFMdT}, but for the sake of completeness we
state it in Proposition \ref{p1.1} below and we ouline its proof
in the appendix.

\begin{proposition}\label{p1.1}
Let $p,\ q>1$, $\delta,\mu  >0$ be such that $(H_1)$ holds. If the system
$$ \quad
\begin{array}{l}
-\Delta_{p} u= |v|^\delta\quad  v>0\quad \mbox{in}\quad \RR^N\cr
-\Delta_{q} v=|u|^{\mu }\quad  u>0\quad \mbox{in}\quad \RR^N
\end{array}\leqno(S_\infty)
$$
has no  radially symmetric solution $(u,v)$ in $(C^1(\mathbb R^N)\cap C^2(\mathbb R^N\setminus\{0\})^2$,
then system $(S_R)$ possesses a  nontrivial solution  for any $R>0$.
\end{proposition}
\medskip

\begin{remark}{\rm We do not know if the converse of this proposition is
true as it is in the case of a single equation or the case of the system with $p=q=2$, see \cite{CFM92}, \cite{Mi93},
\cite{PevdV92}, \cite{V}.}
\end{remark}

\begin{remark}{\rm If $p\ge N$, then $-\Delta_{p}u\ge 0$ and $u\ge 0$ in $\mathbb R^N$ imply $u=Const.$,
see \cite{Ni-Se85}, \cite{Ni-Se85b}, \cite{Ni-Se86}. Hence from $-\Delta_{p}u= 0$ it follows that $v=0$, and from the second
equation it follows that $u=0$. Therefore, if $p\ge N$ or / and $q\ge N$ it follows from Proposition \ref{p1.1} that $(S_R)$
possesses at least one solution $(u,v)$. Hence in Theorem \ref{mainth} we may assume without loss of generality that
$\max\{p,q\}<N$.}
\end{remark}

\begin{remark} {\rm In \cite{CMM93} is has been shown that if
\be\label{*} \max\Bigl\{\alpha-\frac
{N-p}{p-1},\beta-\frac{N-q}{q-1}\Bigr\}\geq 0 \ee where
\be\label{defsigma} \alpha=\frac 1{d}[p(q-1)+\delta q]\quad
\beta=\frac 1{d}[q(p-1)+\mu p] \ee and $(H_1)$, then the
assumptions of Proposition \ref{p1.1} are satisfied. Hence in this
case, that is, when \eqref{*} is satisfied, the existence of
solutions to $(S_R)$ follows. Observe that in case that $p=q=m$
and $\delta =\mu $, the condition \eqref{*} is equivalent to
$$\frac{m-1}{\delta +m-1}>\frac{N-m}{m},\quad m>1$$
which is more restrictive than \eqref{eqcond}. Hence condition \eqref{*} is not optimal.}
\end{remark}

We are now in a position to state our main results.
\begin{theorem}\label{mainth2}
Suppose  $N\ge 2$ and that $\delta $, $\mu >0$ satisfy $(H_1)$.
\begin{enumerate}
\item Let
\be\label{1.2} \frac{2N}{N+1}< p\le 2\mbox{\quad and\quad }\frac{2N}{N+1}< q\le 2.\quad
 \ee
Then  problem $(S_R)$ possesses a solution $(u,v)$ provided that
\be\label{1.3} \frac{1}{\delta +1}+\frac{1}{\mu
+1}>\frac{N-\underline{m}}{N(\underline{m}-1)} \ee where
$\underline{m}=\min\{p,q\}$. \item Let \be\label{1.22} 2\le
p<N\quad\mbox{and}\quad2\le q<N.\quad \ee Then  problem $(S_R)$
possesses a   solution $(u,v)$ provided that \be\label{1.32}
\frac{1}{\delta +1}+\frac{1}{\mu
+1}>\frac{N(\overline{m}-1)-\overline{m}}{N(\overline{m}-1)} \ee
where $\overline{m}=\max\{p,q\}$.
\end{enumerate}
\end{theorem}

\begin{remark} {\rm Observe that when $p=q=2$, condition \eqref{1.3} and \eqref{1.32} are
the same and they are optimal, see \cite{CFM92}, \cite{Mi93}, \cite{PevdV92}, \cite{V}. When $m=p=q\not=2$, $m<2$ and $\delta
=\mu $, condition \eqref{1.3} reads
$$\frac{2}{\delta +1}>\frac{N- m}{N( m-1)}.$$
Since \eqref{eqcond}, which is  optimal, can be rewritten as
$$\frac{2}{\delta +1}>\frac{N- m}{N( m-1)}\frac{2(m-1)}{m},\quad \mbox{with }\ \frac{2(m-1)}{m}<1,$$
it follows that condition \eqref{1.3} is not optimal.

When $m<2$, we note that \eqref{1.3} gives a new region
of existence provided that
$${\frac {N \left( m-1 \right) }{N-m}}<{\frac { \left( 2\,N+1 \right) m-
3\,N}{N-m}},$$ which holds if $m>\frac{2N}{N+1}.$ Since
$2N/(N+1)<2$, there is always room for some $m<2$, as is shown in
Figure \ref{Nm1}.

\begin{figure}[htb]
\begin{center}
\input{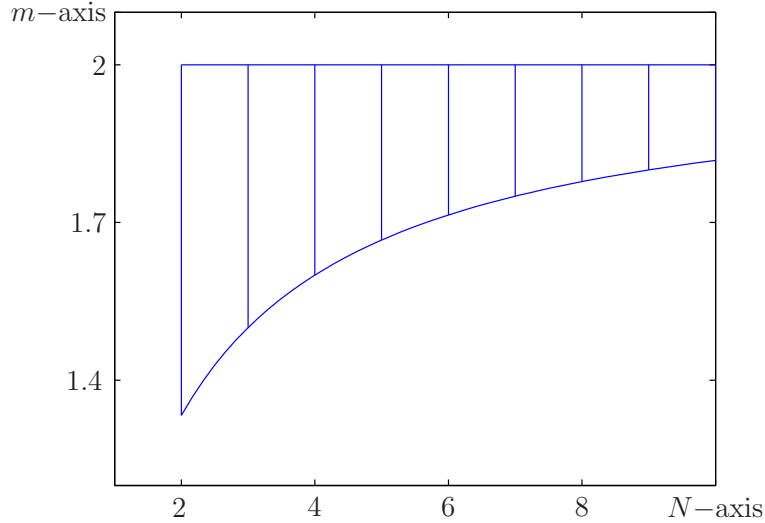}
\end{center}
  \caption{ The values of $m<2$ for which we obtain a new region of existence.}\label{Nm1}
\end{figure}

When $m>2$, we note that \eqref{1.32} gives a new region
of existence provided that
$$\frac{N(m-1)}{N-m}<\frac{N(m-1)+m}{N(m-1)-m},$$
which holds if
$$m<\,{\frac { \left( 3\,N+1+\sqrt {(N-1)(N+7)} \right) N}{2{N}^{2
}+2}}.
$$
Since the right hand side of this inequality is greater than $2$ for $N>2$, there is always room for some $m>2$, as is shown in
the next figure.
\begin{figure}[htb]
\begin{center}
\input{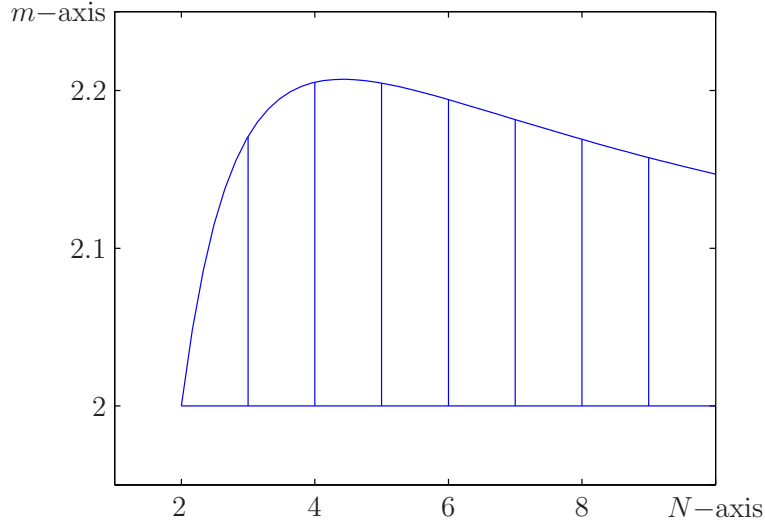}
\end{center}
 \caption{ The values of $m>2$ for which we obtain a new region of
 existence.}\label{Nm2}
\end{figure}
}\end{remark}
\bigskip

\bigskip

Finally we have

\begin{theorem}\label{mainth}
Suppose that $N>2$  and $\delta ,\mu >0$ satisfy  $(H_1)$.
\begin{enumerate}
\item If $2\le p, \ q<N$,  and \be\label{1.4} \frac{1}{\delta
+1}+\frac{1}{\mu +1}\le\frac{N- \overline{m}}{N( \overline{m}-1)},
\ee where $\overline{m}=\max\{p,q\}$, then  system $(S_R)$ has no
solutions (regular or not). \item If $N/(N-1)<p,\ q\le 2$ and
\be\label{1.42} \frac{1}{\delta +1}+\frac{1}{\mu
+1}\le\frac{N(\underline{m}-1)-\underline{m}}{N(\underline{m}-1)},
\ee where $\underline{m}=\min\{p,q\}$, then  system $(S_R)$ has no
  solutions (regular or not).
\end{enumerate}
\end{theorem}

\begin{remark} {\rm In case that $p=q=2$, \eqref{1.4}-\eqref{1.42} are optimal, but when $m=p=q\not=2$ they are not.
Indeed,
if we set $\mu=\delta$ in \eqref{1.4}, we obtain
$$\mu\ge \frac{(2N+1)m-3N}{N-m}>\frac{N(m-1)+m}{N-m}\quad\mbox{for all }\ m>2,$$
and if we set $\mu=\delta$ in \eqref{1.42}, we obtain
$$\mu\ge \frac{N(m-1)+m}{N(m-1)-m}>\frac{N(m-1)+m}{N-m}\quad\mbox{for all }\ m<2.$$
Since $\mu\ge \frac{N(m-1)+m}{N-m}$ is the optimal range for the
case of one equation, our claim follows.}
\end{remark}

\begin{figure}[htb]
\begin{center}
%
%
\begin{psfrags}%
\psfragscanon%
%
\psfrag{s04}[l][l]{\color[rgb]{0,0,0}\setlength{\tabcolsep}{0pt}\begin{tabular}{l}{\large NE}\end{tabular}}%
\psfrag{s05}[l][l]{\color[rgb]{0,0,0}\setlength{\tabcolsep}{0pt}\begin{tabular}{l}{\large E}\end{tabular}}%
\psfrag{s06}[b][b]{\color[rgb]{0,0,0}\setlength{\tabcolsep}{0pt}\begin{tabular}{c}$m=1.9$ and $N=4$\end{tabular}}%
%
\psfrag{x01}[t][t]{0}%
\psfrag{x02}[t][t]{}%
\psfrag{x03}[t][t]{}%
\psfrag{x04}[t][t]{$\frac{N(m-1)+m}{N-m}$}%
\psfrag{x05}[t][t]{$\mu-$axis}%
%
\psfrag{v01}[r][r]{$\frac{m(m-1)}{N-m}$}%
\psfrag{v02}[r][r]{$\frac{N(m-1)}{N-m}$}%
\psfrag{v03}[r][r]{}%
\psfrag{v04}[r][r]{$\delta-$axis}%
%
\resizebox{9cm}{!}{\includegraphics{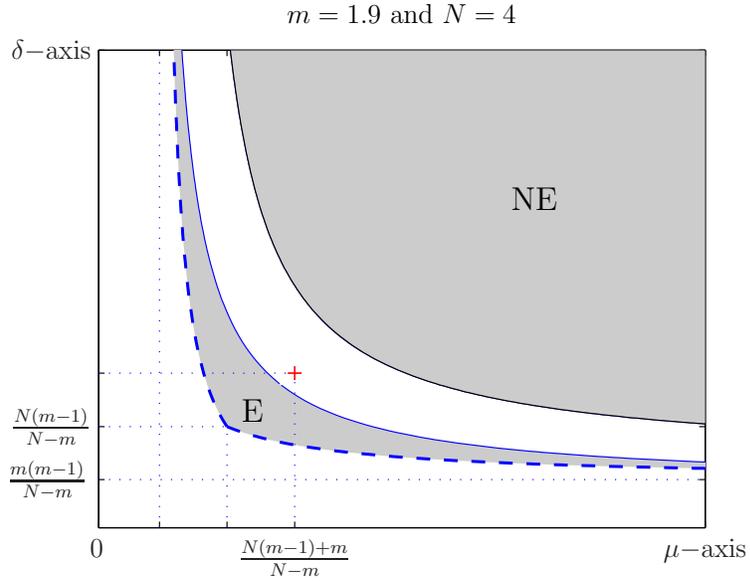}}%
\end{psfrags}%
%

\end{center}
\caption{The nonexistence region in grey is given by \eqref{1.42},
and the new existence  region is the one bounded by equality in
\eqref{*} and by equality in \eqref{1.3}.} \label{intro2}
\end{figure}

\begin{figure}[htb]
\begin{center}
%
%
\begin{psfrags}%
\psfragscanon%
%
\psfrag{s04}[l][l]{\color[rgb]{0,0,0}\setlength{\tabcolsep}{0pt}\begin{tabular}{l}{\large NE}\end{tabular}}%
\psfrag{s05}[l][l]{\color[rgb]{0,0,0}\setlength{\tabcolsep}{0pt}\begin{tabular}{l}{\large E}\end{tabular}}%
\psfrag{s06}[b][b]{\color[rgb]{0,0,0}\setlength{\tabcolsep}{0pt}\begin{tabular}{c}$m=2.1$ and $N=4$\end{tabular}}%
%
\psfrag{x01}[t][t]{0}%
\psfrag{x02}[t][t]{}%
\psfrag{x03}[t][t]{}%
\psfrag{x04}[t][t]{$\frac{N(m-1)+m}{N-m}$}%
\psfrag{x05}[t][t]{$\mu-$axis}%
%
\psfrag{v01}[r][r]{$\frac{m(m-1)}{N-m}$}%
\psfrag{v02}[r][r]{$\frac{N(m-1)}{N-m}$}%
\psfrag{v03}[r][r]{}%
\psfrag{v04}[r][r]{$\delta-$axis}%
%
\resizebox{9cm}{!}{\includegraphics{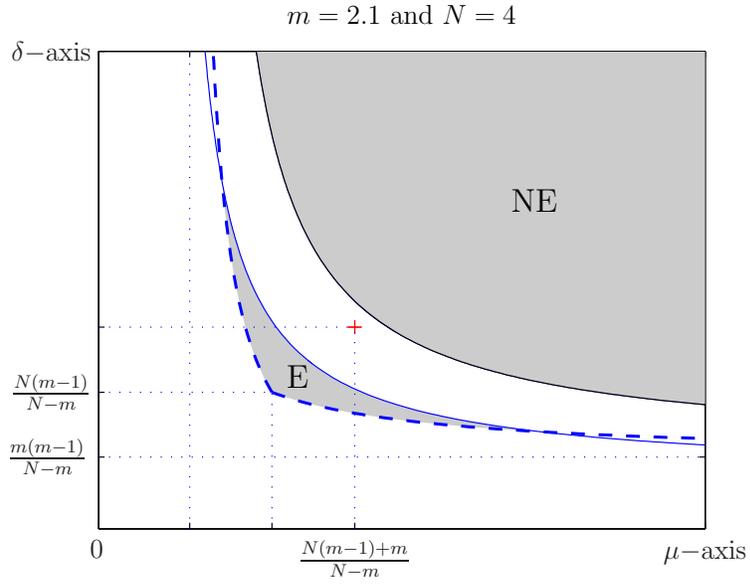}}%
\end{psfrags}%
%

\end{center}
\caption{  The non existence region in grey is given by
\eqref{1.4}. The new existence region is bounded  by equality in
\eqref{1.32} and the dashed curve which is given by equality
\eqref{*}.} \label{intro1}
\end{figure}

For $N=4$, in Figure \ref{intro2} we show the new region of
existence  and the non-existence region for the case $m=p=q=1.9$,
and in Figure \ref{intro1} we show the new region of existence and
the non-existence region for the case $m=p=q=2.1$.

\clearpage
\newpage

Our article is organized as follows. In section \ref{pohozaev} we give a Pohozaev type identity which
is the key to prove our
main results and we prove them.
Finally in the appendix we
prove Proposition \ref{p1.1}.

\section{A Pohozaev type identity and proof of our main results}\label{pohozaev}

Our main theorems are based on the two following lemmas, which
give appropriate generalizations of the Pohozaev identity used to
deal with the case $p=q=2$, see \cite{Mi93}.
\medskip

\begin{lemma}\label{Poh1} Let $(u,v)\in (C^1[0,\infty)\cap C^2(0,\infty))^2 $ be a  solution of the system
$$ \quad
\begin{array}{l}
-(r^{N-1}|u'|^{p-2}u' )'= r^{N-1}v^\delta  \cr
-(r^{N-1}|v'|^{q-2}v' )'=r^{N-1}u^\mu\cr
u(r)>0,\quad v(r)>0,\quad r\in[0,\infty),
\end{array}\leqno(S^r_{\infty})
$$
with $\delta,\ \mu>0$, and assume that either
$$2N/(N+1)\le p\le q\le 2\quad{ or }\quad
2\le p\le q<N.$$

Let us define
\begin{eqnarray}
E_1(r)=r^{N+k_1-2}|u'|^{p-1}|v'|^{q-1}\qquad\qquad\qquad\qquad\qquad\qquad \nonumber\\
-\frac{N}{\delta +1}r^{N-1}|u'|^{p-1}\int_r^{\infty}s^{k_1-2}|v'|^{q-1}\,ds
-\frac{N}{\mu +1}r^{N-1}|v'|^{q-1}\int_r^{\infty}s^{k_1-2}|u'|^{p-1}\,ds \nonumber\\
+r^{N}\int_r^{\infty}s^{k_1-2}|v'|^{q-1}v^\delta \,ds+
r^N\int_r^{\infty}s^{k_1-2}|u'|^{p-1}u^\mu \,ds,\, r\in
(0,\infty),\quad\label{energy1}
\end{eqnarray}
where
$$k_1=\begin{cases}p+\frac{N-p}{p-1}(p-2)\quad\mbox{ if }2N/(N+1)\le p\le q\le 2,\\
\frac{q}{q-1}\quad\mbox{ if } 2\le p\le q.
\end{cases}$$
Then for $r\in (0,\infty)$ we have
\begin{eqnarray}\label{Ep1}
E_1'(r)&=&\Bigl(k_1-N +\frac{N}{\delta +1}+\frac{N}{\mu +1}\Bigr)r^{N+k_1-3}|u'|^{p-1}|v'|^{q-1}\nonumber \\
 &&-\frac{N}{\delta +1}r^{N-1}v^{\delta }\int_r^{\infty}s^{k_1-2}|v'|^{q-1}\,ds+
Nr^{N-1}\int_r^{\infty}s^{k_1-2}|v'|^{q-1}v^\delta \,ds \nonumber\\
 &&-\frac{N}{\mu +1}r^{N-1}u^{\mu }\int_r^{\infty}s^{k_1-2}|u'|^{p-1}\,ds+
 Nr^{N-1}\int_r^{\infty}s^{k_1-2}|u'|^{p-1}u^\mu \,ds.\qquad
\end{eqnarray}
\end{lemma}

\begin{proof} We prove first that $E_1$ is well
defined.
Since $u$ and $v$ are decreasing functions we have
that for any $T>r$ it holds that
$$\int_r^{T}s^{k_1-2}|v'|^{q-1}v^\delta \,ds\le
v^\delta (r)\int_r^{T}s^{k_1-2}|v'|^{q-1}\,ds,$$ and
$$
\int_r^{T}s^{k_1-2}|u'|^{p-1}u^\mu \,ds\le
u^\mu (r)\int_r^{T}s^{k_1-2}|u'|^{p-1}ds,$$   thus it is sufficient to prove that the first two integrals
in \eqref{energy1} are well defined.

In order to do so, we recall that from  \cite[Lemma 2.1]{CMM93} or \cite[Proposition V.1]{CFMdT}, we have
$$ |u'(r)|\leq C r^{-\alpha-1}\quad
|v'(r)|\leq C r^{-\beta-1}
$$
for some $C>0$ and $r$ large.

We deal first with the case   $2N/(N+1)\le p\le q\le 2$. We will
see first that $1-k_1+(\alpha+1)(p-1)>0$, which is the condition
to have the second integral in \eqref{energy1} well defined.
Indeed,
$$1-k_1+(\alpha+1)(p-1)=p-k_1+\alpha(p-1)=\frac{N-p}{p-1}(2-p)+\alpha(p-1)>0.$$
For the first integral we have
$$1-k_1+(\beta+1)(q-1)=q-k_1+\beta(q-1)\ge p-k_1+\beta(q-1) >0,$$
 hence the first integral appearing in \eqref{energy1} is also well defined.

 For the case $2\le p\le q$, it can be easily verified that
$\alpha(p-1)>k_1-1-(p-1)$ and $\beta(q-1)>k_1-1-(q-1)$.
We only verify the first inequality: as $p\ge 2$,  we have that
$$k_1-1-(p-1)=\frac{1-(q-1)(p-1)}{q-1}\le 0,$$
and $\alpha(p-1)>0$, thus the first two integrals in \eqref{energy1} are well defined.

Now \eqref{Ep1} follows  by direct differentiation using that
$(u,v)$ is a   solution to $(S^r_{\infty})$.
\end{proof}

Our second lemma is essentially the same Pohozaev identity as in
Lemma \ref{Poh1}, but in $(0,R]$.
\begin{lemma}\label{Poh2} Let $(u,v)\in (C^1[0,R]\cap C^2(0,R])^2$ be a  solution of the
system $(S_R)$ with $\delta,\ \mu>0$,
 and assume that either
 $$N/(N-1)< p\le q\le 2\quad\mbox{ or  }\quad 2\le p\le q<N.$$
Let us define
\begin{eqnarray}
E_2(r)=r^{N+k_2-2}|u'|^{p-1}|v'|^{q-1}\qquad\qquad\qquad\qquad\qquad\qquad \nonumber\\
-\frac{N}{\delta +1}r^{N-1}|u'|^{p-1}\int_r^{R}s^{k_2-2}|v'|^{q-1}\,ds
-\frac{N}{\mu +1}r^{N-1}|v'|^{q-1}\int_r^{R}s^{k_2-2}|u'|^{p-1}\,ds \nonumber\\
+r^{N}\int_r^{R}s^{k_2-2}|v'|^{q-1}v^\delta \,ds+
r^N\int_r^{R}s^{k_2-2}|u'|^{p-1}u^\mu \,ds,\, r\in
(0,R],\quad\label{energy2}
\end{eqnarray}
where
$$k_2=\begin{cases}q+\frac{N-q}{q-1}(q-2)\quad\mbox{ if }2\le p\le q<N, \\
\frac{p}{p-1}\quad\mbox{ if } N/(N-1)< p\le q\le 2.
\end{cases}$$
Then for $r\in (0,R)$ we have
\begin{eqnarray}\label{Ep2}
E_2'(r)&=&\Bigl(k_2-N +\frac{N}{\delta +1}+\frac{N}{\mu +1}\Bigr)r^{N+k_2-3}|u'|^{p-1}|v'|^{q-1}\nonumber \\
 &&-\frac{N}{\delta +1}r^{N-1}v^{\delta }\int_r^{R}s^{k_2-2}|v'|^{q-1}\,ds+
Nr^{N-1}\int_r^{R}s^{k_2-2}|v'|^{q-1}v^\delta \,ds \nonumber\\
 &&-\frac{N}{\mu +1}r^{N-1}u^{\mu }\int_r^{R}s^{k_2-2}|u'|^{p-1}\,ds+
 Nr^{N-1}\int_r^{R}s^{k_2-2}|u'|^{p-1}u^\mu \,ds.\qquad
\end{eqnarray}
\end{lemma}

\medskip

Now we can prove our main results.

\begin{proof}[Proof of Theorem \ref{mainth2}]
In view of Proposition \ref{p1.1}, in order to prove our theorem
we only need to prove that under assumption \eqref{1.3} or
\eqref{1.32} system $(S_\infty)$ does not possess any   radial
solution. We will argue by contradiction by assuming that there
exists a   radially symmetric solution $(u,v)$ to $(S_\infty)$.
The idea is to have $E_1$ strictly increasing with
$\lim\limits_{r\to 0^+}E_1(r)=0$ and $\lim\limits_{r\to
\infty}E_1(r)=0$ which will give a contradiction.

We prove first (1), and start by proving
that $\lim\limits_{r\to 0}E_1(r)=0.$ Since $u$ and $v$ are regular, a
simple application of L'H\^opital's rule gives
$$
\lim\limits_{r\to
0}|u'(r)|^{p-1}/r=v(0)^\delta /N\quad\mbox{and}\quad\lim\limits_{r\to
0}|v'(r)|^{q-1}/r=u(0)^\mu /N.
$$
Therefore we need $ N+k_1>0$, or equivalently,  $p>3N/(2N+1)$.
But $p>2N/(N+1)>3N/(2N+1)$ if $N>1$, hence $\lim\limits_{r\to
0}E_1(r)=0$.

 We now verify that $\lim\limits_{r\to \infty}E_1(r)=0.$ We have the
bounds near infinity given by
$$
u(r)\le Cr^{-\alpha},\quad |u'(r)|\leq C r^{-\alpha-1},\quad
v(r)\le Cr^{-\beta},\quad |v'(r)|\leq C r^{-\beta-1}
$$
for some $C>0$ and $r$ large, see \cite[Lemma 2.1]{CMM93} or
\cite[Proposition V.1]{CFMdT}. Next, by observing that by the
definition of $\alpha,\ \beta$ we have
$$1-\delta \beta=-(\alpha+1)(p-1)\quad\mbox{and}\quad
1-\mu \alpha=-(\beta+1)(q-1),$$ in order  that
$\lim\limits_{r\to \infty}E(r)=0$ it is sufficient to show that
$$
N+k_1-2-(\alpha+1)(p-1)-(\beta+1)(q-1)<0.
$$

This last inequality is equivalent to
$$
N+k_1-p-q< \frac{(p-1)[p(q-1)+\delta q]
+(q-1)[q(p-1)+\mu p]}{\delta \mu -(p-1)(q-1)}.
$$
Calling $A=\delta +q-1,$  $B=\mu +p-1$ and $L=N+k_1-p-q,$ this
reads
$$
L<\frac{q(p-1)A+p(q-1)B}{AB-(p-1)A-(q-1)B},
$$
and since the denominator is positive, we have then to prove that
\be\label{toprove} \quad
L<\frac{(q-1)(L+p)}{\delta +q-1}+\frac{(p-1)(L+q)}{\mu +p-1},
\ee
 Since $p\le q\le 2$, we have that $\delta +q-1\le \delta +1$ and $\mu +p-1\le
 \mu +1$, and thus, using assumption \eqref{1.3} (with
 $\underline{m}=p$),
we have (using also that $0<L+p<L+q$)
$$
\frac{(L+p)(N-p)}{N}\leq\frac{(q-1)(L+p)}{\delta +q-1}+\frac{(p-1)(L+q)}{\mu +p-1}.
$$
Therefore we have  to prove that $L<(L+p)(N-p)/N$, which is
equivalent to $k_1-q<0$. Using now that $p<2$, we have
\begin{eqnarray*}
 k_1=p+
\frac{N-p}{p-1}(p-2)<q,
\end{eqnarray*}
proving \eqref{toprove} and thus $E_1(\infty)=0$.

We prove next that under the assumptions of the theorem we have
$E_1'(r)>0$ for all $r>0$. Since by the choice of $k_1$
$$k_1-N+\frac{N}{\delta +1}+\frac{N}{\mu +1}=-\frac{(N-p)}{p-1}+\frac{N}{\delta +1}+\frac{N}{\mu +1},$$
we have by assumption \eqref{1.3} that the first term in \eqref{Ep1} is indeed positive.
Let us set now
\begin{equation}\label{pr1g}
G(p,\mu ,u)(r)=N\int_r^\infty
s^{k_1-2}|u'|^{p-1}u^\mu \,ds-\frac{N}{\mu +1}u^{\mu }\int_r^\infty
s^{k_1-2}|u'|^{p-1}\,ds,
\end{equation}
where $k_1=p+
\frac{N-p}{p-1}(p-2)$. With this notation, for $r\in (0,\infty),$ we have that $E_1'(r)$ can be written as
\begin{eqnarray}\label{EEp1}
E_1'(r)&=&\Bigl(k_1-N +\frac{N}{\delta +1}+\frac{N}{\mu +1}\Bigr)r^{N+k_1-3}|u'|^{p-1}|v'|^{q-1}\nonumber \\
 &&+r^{N-1}G(p,\mu ,u)(r)+r^{N-1}G(q,\delta ,v)(r).\qquad
\end{eqnarray}

By differentiating both sides in \eqref{pr1g} with respect to $r$ we obtain
\begin{equation}\label{pr21}
G'(p,\mu ,u)(r)=\Bigl(-N+\frac{N}{\mu +1}\Bigr)r^{k_1-2}|u'|^{p-1}u^\mu +
\frac{N\mu }{\mu +1}u^{\mu -1}|u'|\int_r^\infty
s^{k_1-2}|u'|^{p-1}\,ds.
\end{equation}
Using now that $(r^{N-1}|u'|^{p-1})'\geq 0$, we have
$s^{N-1}|u'|^{p-1}(s)\geq  r^{N-1}|u'|^{p-1}(r)$ for $s\geq
r$, and consequently, using that  $p\leq 2,$ we find that
$\bigl(s^{(N-1)/(p-1)}|u'|(s)\bigr)^{p-2}\leq
\bigl(r^{(N-1)/(p-1)}|u'|(r)\bigr)^{p-2}$ for $s\geq r.$
Therefore,
\begin{equation}\label{pr31}
\int_r^\infty
s^{k_1-2}|u'|^{p-1}\,ds=\int_r^\infty\bigl(s^{(N-1)/(p-1)}|u'|(s)\bigr)^{p-2}|u'(s)|\,ds\leq
r^{k_1-2}u(r)|u'|^{p-2}.
\end{equation}
Replacing (\ref{pr31}) into (\ref{pr21}), we obtain
$$
G'(p,\mu ,u)(r)\leq
\Bigl(-N+\frac{N}{\mu +1}(\mu +1)\Bigr)r^{k-2}|u'|^{p-1}u^\mu =0,
$$
hence $G'(p,\mu ,u)(r)\leq 0$ for all $r>0$, and since
$G(p,\mu ,u)(\infty)=0,$ we have \linebreak[5] $G(p,\mu ,u)(r)\geq 0$ for all
$r>0$. Thus the  term in the third line in \eqref{Ep1} is also
positive.

Finally we show that $G(q,\delta ,v)(r)\geq 0$ for all $r>0$, proving
that the  term in the second line of \eqref{Ep1} is also positive.
Indeed, we define $\bar k=q+\frac{N-q}{q-1}(q-2)$ and note
that $k_1\leq \bar k$ when $q\geq p.$ We proceed as above
using the following inequality
$$
\int_r^\infty s^{k_1-2}|v'|^{q-1}\,ds=\int_r^\infty s^{k_1-\bar
k}\bigl(s^{(N-1)/(q-1)}|v'|(s)\bigr)^{q-2}|v'(s)|\,ds\leq
$$
$$
r^{k_1-\bar k}\int_r^\infty
\bigl(s^{(N-1)/(q-1)}|v'|(s)\bigr)^{q-2}|v'(s)|\,ds\leq
r^{k_1-2}v(r)|v'|^{q-2}.
$$

Therefore $E_1'(r)>0$ for all $r>0$ in contradiction with
$E_1(0^+)=E_1(\infty)=0$. Thus under the assumptions of the
Theorem there cannot exist    radially symmetric solutions to
$(S_\infty)$ and we can use Proposition \ref{p1.1} to obtain the
existence of at least one   solution to $(S_R)$ for any positive
$R$.
\medskip

We next prove (2) and hence we assume $q=\overline{m}$.
The proof of $E_1(0^+)=0$ follows as before. In order to prove that $E_1(\infty)=0$ we need to prove
\eqref{toprove} in the case $k_1=q/(q-1)$. Since $p,\ q\ge 2$, and the function
$x\mapsto x/(c+x)$ is strictly increasing in $(0,\infty)$ for any $c>0$, we have that
$$\frac{p-1}{\mu +p-1}\ge \frac{1}{\mu +1}\quad\mbox{and}\quad
\frac{q-1}{\delta +q-1}\ge \frac{1}{\delta +1}.$$
Since $q<N$, we have that $L+p=N+k_1-q>0$, hence by assumption \eqref{1.32}, we have
\begin{eqnarray*}
\frac{(q-1)(L+p)}{\delta +q-1}+\frac{(p-1)(L+q)}{\mu +p-1}&\ge&
(L+p)\Bigl(\frac{q-1}{\delta +q-1}+\frac{p-1}{\mu +p-1}\Bigr)\\
&\ge&(L+p)\Bigl(\frac{1}{\delta +1}+\frac{1}{\mu +1}\Bigr)\\
&>& (L+p)\Bigl(1-\frac{k_1}{N}\Bigr),
\end{eqnarray*}
and therefore \eqref{toprove} will follow if we prove that
\begin{equation}\label{last-tp}
(L+p)\Bigl(1-\frac{k_1}{N}\Bigr)\ge L.
\end{equation}
Now, \eqref{last-tp} is equivalent to $(N+k_1-q)k_1\le Np$. Since $k_1\le 2$, $p\ge 2$,
and $N+k_1-q\le N$, \eqref{last-tp} follows and $E_1(\infty)=0$.
\medskip

We prove next that $E_1'(r)>0$ for all $r>0$.
Now the first term in \eqref{Ep1} is positive by assumption \eqref{1.32}.
We set as before
\begin{equation}\label{pr1g2}
 G(p,\mu ,u)(r)=N\int_r^\infty
s^{k_1-2}|u'|^{p-1}u^\mu \,ds-\frac{N}{\mu +1}u^{\mu }\int_r^\infty
s^{k_1-2}|u'|^{p-1}\,ds,
\end{equation}
where now $k_1= {q}/(q-1)$, obtaining again that
\begin{eqnarray}\label{pr212}
 G'(p,\mu ,u)(r)=\Bigl(-N+\frac{N}{\mu +1}\Bigr)r^{k_1-2}|u'|^{p-1}u^\mu \\
+\frac{N\mu }{\mu +1}u^{\mu -1}|u'|\int_r^\infty
s^{k_1-2}|u'|^{p-1}\,ds.
\end{eqnarray}
We claim that $|u'|^{p-1}/r$ is decreasing for all $r>0$: indeed, since
$$\frac{|u'|^{p-1}}{r}=\frac{1}{r^N}\int_0^r s^{N-1}v^\delta (s)ds,$$
we have that
$$
\frac{d}{dr}\frac{(|u'|^{p-1})}{r}=\frac{1}{r^N} r^{N-1}v^\delta (r)-N\frac{r^{N-1}}{r^{2N}}
\int_0^r s^{N-1}v^\delta (s)ds,
$$
and thus, using that $v$ is decreasing in $(0,\infty)$ we find that
\begin{equation}\label{dec-quot}
\frac{d}{dr}\frac{(|u'|^{p-1})}{r}\le \frac{1}{r^N} r^{N-1}v^\delta (r)-N\frac{r^{N-1}}{r^{2N}}
\frac{r^{N}}{N}v^\delta (r)=0.
\end{equation}
Since for $\bar k=p/(p-1)$, it holds that
$$s^{k_1-2}|u'|^{p-1}=s^{k_1-\bar k}\Bigl(\frac{|u'|}{s^{1/(p-1)}}\Bigr)^{p-2}|u'|,$$
and since $p\le q$, also $k_1\le \bar k$. Hence we find that
$$\int_r^\infty
s^{k_1-2}|u'|^{p-1}\,ds\le r^{k_1-2}|u'|^{p-2}\int_r^\infty |u'(s)|ds=r^{k_1-2}|u'|^{p-2}u(r).$$
Thus, replacing this estimate into  \eqref{pr212} we obtain that
$$ G'(p,\mu ,u)(r)\le \Bigl(-N+\frac{N(\mu +1)}{\mu +1}\Bigr)r^{k_1-2}|u'|^{p-1}u^\mu = 0.
$$
Since
$ G(p,\mu ,u)(\infty)=0,$ we have $ G(p,\mu ,u)(r)\geq 0$ for all
$r>0$. Thus the  term in the third line of \eqref{Ep1} is
positive.

The same argument, with $\bar k=k_1$ can be used to show that $ G(q,\delta ,v)(r)\geq 0$ for all $r>0$, proving
that the  term in the second line of \eqref{Ep1} is also positive and thus $E_1'(r)>0$ for
all $r>0$. Again we obtain a contradiction and we can use Proposition \ref{p1.1} to obtain the existence
of at least one
  solution to $(S_R)$ for any positive $R$.

\end{proof}

Finally in this section we prove Theorem \ref{mainth}.

\begin{proof}[Proof of Theorem \ref{mainth}]We will argue by contradiction assuming that there exists
a   solution $(u,v)$ to $(S_R)$. Now we will use Lemma \ref{Poh2}.
The idea is to have $E_2$ decreasing with $E_2(0^+)=0$ and
$E_2(R)>0$ yielding a contradiction.

We assume $q=\overline{m},$ $p=\underline{m},$ and
since the case $p=q=2$ was proven
in \cite{Mi93}, we may assume without loss of generality that $q>2$ for part (1) and $p<2$ for part (2).

By direct computation we have that
$$E_2(R)=R^{N+k_2-2}|u'(R)|^{p-1}|v'(R)|^{q-1}>0.$$

We will show next that $E_2(0^+)=0$.
By  \cite[Lemma 2.1]{CMM93} or
 \cite[Proposition V.1]{CFMdT},  we have that any   solution $(u,v)$ to $(S_R)$ satisfies
$$
u(r)\le K r^{-\alpha},\quad |u'(r)|\leq K r^{-\alpha-1},\quad
v(r)\le K r^{-\beta},\quad |v'(r)|\leq K r^{-\beta-1}
$$
for some $K>0$ and $0<r\ll 1$.  Hence in order to show that $E_2(0^+)=0$ we need
\be\label{toprove2}
N+k_2-2-(\alpha+1)(p-1)-(\beta+1)(q-1)>0.
\ee
 As for \eqref{toprove}, this last inequality
reduces to
\begin{equation*}
 L>\frac{(q-1)(L+p)}{\delta +q-1}+\frac{(p-1)(L+q)}{\mu +p-1}.
\end{equation*}
where $L:=N+k_2-p-q.$

We deal first with the case $2\le p\le q$. Since in this
case $\overline{m}= q$, by assumption
\eqref{1.4} we have
$$
\frac{N-q}{(q-1)N}\geq\frac{1}{\delta +1}+\frac{1}{\mu +1}.
$$
On the other hand, using that $\delta +q-1\ge \delta +1$ and $\mu +p-1\ge
\mu +1$, and $q\ge p$ we have
$$
(q-1)(L+q)\Bigl(\frac{1}{\delta +1}+\frac{1}{\mu +1}\Bigr)\ge\frac{(q-1)(L+p)}{\delta +q-1}+\frac{(p-1)(L+q)}{\mu +p-1},
$$
which implies
$$
\frac{(L+q)(N-q)}{N}\geq\frac{(q-1)(L+p)}{\delta +q-1}+\frac{(p-1)(L+q)}{\mu +p-1}.
$$
Hence in order to prove \eqref{toprove2} it is sufficient that $L>(L+q)(N-q)/N$. But this is
equivalent to prove that $q+
\frac{N-q}{q-1}(q-2)>p$, which is clearly true by the assumption $q> 2$ and $q\ge p$.
\medskip

Next we deal with the case $N/(N-1)<p\le q\le 2$. Using again the monotonicity of the function
$x\mapsto x/(c+x)$, where $c>0$, using now that $p-1\le 1$ and $q-1\le 1$, we find that
$$\frac{p-1}{\mu +p-1}\le \frac{1}{\mu +1}\quad\mbox{and}\quad
\frac{q-1}{\delta +q-1}\le \frac{1}{\delta +1}.$$
Hence by assumption \eqref{1.42} and using that $L+q>0$ we find that
\begin{eqnarray*}
\frac{(q-1)(L+p)}{\delta +q-1}+\frac{(p-1)(L+q)}{\mu +p-1}&\le&
(L+q)\Bigl(\frac{q-1}{\delta +q-1}+\frac{p-1}{\mu +p-1}\Bigr)\\
&\le&(L+q)\Bigl(\frac{1}{\delta +1}+\frac{1}{\mu +1}\Bigr)\\
&\le& (L+q)\Bigl(1-\frac{k_2}{N}\Bigr).
\end{eqnarray*}
Hence in order to establish \eqref{toprove2}, it is sufficient that
$$(L+q)\Bigl(1-\frac{k_2}{N}\Bigr)<L.$$
Since this inequality is equivalent to $Nq<(N+k_2-p)k_2$,
and $q<N$, $k_2\ge 2$, a sufficient condition so that it holds is that
$N<N+k_2-p$, which is clearly true since $p<2$ and $k_2=p/(p-1)$.
\medskip

Finally we prove that $E_2'(r)\le 0$ for all $r\in(0,R)$.
To this end  we define
\begin{equation}\label{pr1}
G(q,\delta ,v)(r)=N\int_r^Rs^{k_2-2}|v'|^{q-1}v^\delta \,ds-\frac{N}{\delta +1}v^{\delta }
\int_r^Rs^{k_2-2}|v'|^{q-1}\,ds,
\end{equation}
so  that
\begin{eqnarray}\label{pr1a}
E_2'(r)=r^{N-1}\Bigl(k_2-N
+\frac{N}{\delta +1}+\frac{N}{\mu +1}\Bigr)r^{k_2-2}|u'|^{p-1}|v'|^{q-1}\nonumber\\
+r^{N-1}\Bigl(G(q,\delta ,v)(r)+G(p,\mu ,u)(r)\Bigr).\qquad\qquad
\end{eqnarray}
We claim  that $E_2'(r)\leq 0$ for $r\in (0,R)$. Indeed,
differentiating in (\ref{pr1}) with respect to $r$, we obtain
\begin{equation}\label{pr2}
G'(q,\delta ,v)(r)=\Bigl(-N+\frac{N}{\delta +1}\Bigr)r^{k_2-2}|v'|^{q-1}v^\delta +\delta
\frac{N}{\delta +1}v^{\delta -1}|v'|
\int_r^Rs^{k_2-2}|v'|^{q-1}\,ds.
\end{equation}
Assume first that $2\le p\le q$. Using that $(r^{N-1}|v'|^{q-1})'\geq 0$ we have
$$s^{N-1}|v'|^{q-1}(s)\geq  r^{N-1}|v'|^{q-1}(r)\quad\mbox{ for $s\geq
r$},$$
 and consequently,
$\bigl(s^{(N-1)/(q-1)}|v'|(s)\bigr)^{q-2}\geq
\bigl(r^{(N-1)/(q-1)}|v'|(r)\bigr)^{q-2}$ for $s\geq r.$ Hence
using that $k_2-2=\frac{N-1}{q-1}(q-1)$ we obtain
\begin{equation}\label{pr3}
\int_r^Rs^{k_2-2}|v'|^{q-1}\,ds=\int_r^R\bigl(s^{(N-1)/(q-1)}|v'|(s)\bigr)^{q-2}|v'(s)|\,ds\geq
r^{k_2-2}v(r)|v'|^{q-2}.
\end{equation}
Thus replacing (\ref{pr3}) into (\ref{pr2}), we get
$$
G'(q,\delta ,v)(r)\geq
\Bigl(-N+\frac{N}{\delta +1}(\delta +1)\Bigr)r^{k_2-2}|v'|^{q-1}v^\delta =0,
$$
implying $G'(q,\delta ,v)(r)\geq 0$ for all $r\in(0,R)$.

Similarly we obtain that $G'(p,\mu ,u)(r)\geq 0$ for all
$r\in(0,R)$. Indeed, we define $\bar
k=p+\frac{N-p}{p-1}(p-2)$ and note that $k_2\geq \bar k$
when $q\geq p.$ We proceed as before, but using the
inequality
$$
\int_r^Rs^{k_2-2}|u'|^{p-1}\,ds=\int_r^Rs^{k_2-\bar
k}\bigl(s^{(N-1)/(p-1)}|u'|(s)\bigr)^{p-2}|u'(s)|\,ds\geq
$$ $$
r^{k_2-\bar
k}\int_r^R\bigl(s^{(N-1)/(p-1)}|u'|(s)\bigr)^{p-2}|u'(s)|\,ds\geq
r^{k_2-2}u(r)|u'|^{p-2}.
$$
This implies $G'(p,\mu ,u)(r)\geq 0$ for all $r\in(0,R)$.

For the case $N/(N-1)< p\le q\le 2$, we argue as follows: we set $\bar k=q/(q-1)\le k_2$ to obtain
\begin{equation}\label{pr32}
\int_r^Rs^{k_2-2}|v'|^{q-1}\,ds=\int_r^Rs^{k_2-\bar k}\Bigl(\frac{|v'|(s)}{s^{1/(q-1)}}\Bigr)^{q-2}|v'(s)|\,ds\geq
r^{k_2-2}v(r)|v'|^{q-2},
\end{equation}
hence in this case we find that
$$G'(q,\delta ,v)(r)\ge \Bigl(-N+\frac{N}{\delta +1}(\delta +1)\Bigr)r^{k_2-2}|v'|^{q-1}v^\delta = 0.
$$
Similarly, setting $\bar k=k_2$,  we find that
$$G'(p,\mu ,u)(r)\ge \Bigl(-N+\frac{N}{\mu +1}(\mu +1)\Bigr)r^{k_2-2}|v'|^{q-1}v^\delta = 0.
$$
Since
$G(q,\delta ,v)(R)=G(p,\mu ,u)(R)=0,$ we conclude $G(q,\delta ,v)(r)\leq 0$ and \linebreak[5]
$G(p,\mu ,u)(r)\le 0$ for all
$r\in(0,R)$.

Now using
that $(\delta ,\mu )$ satisfies \eqref{1.4}, respectively \eqref{1.42}, $G(q,\delta ,v)(r)\leq 0$ and
$G(p,\mu ,u)(r)\leq 0$ for all $r\in(0,R]$, by   (\ref{pr1a}) we obtain $E_2'(r)< 0$ for all $r\in(0,R]$,
which is a contradiction.
\medskip

Thus the theorem follows.
\end{proof}

\section{Appendix}

We start this section by proving Proposition \ref{p1.1}.

\begin{proof}[Proof of Proposition \ref{p1.1}]

In order to prove Proposition \ref{p1.1}, we will make use of the
ideas first used in \cite{CMM93} and  later in \cite{CFMdT}. For
the convenience of the reader, we summarize below the results that
we shall use. To this end, we define the operator $T$ associated
to system $(S_R)$:

\noindent For $(u,v)\in C[0,R]\times C[0,R]$, we set
\begin{eqnarray*}
T(u,v)(r)&=&\Bigl(\int_r^R\Bigl(s^{1-N}\int_0^s t^{N-1}|v(t)|^\delta dt\Bigr)^{1/(p-1)}ds,\\
&&\qquad\quad\int_r^R \Bigl(s^{1-N}\int_0^s t^{N-1}|u(t)|^\mu
dt\Bigr)^{1/(q-1)}ds\Bigr),
\end{eqnarray*}
and denote by $B(0,s)$, $s>0$, the open ball in $C[0,R]\times C[0,R]$ of radius $s$ centered at the origin.

$T$ has the following properties:
\begin{itemize}
\item[(A)]
\begin{enumerate}
\item[(i)] $T$ maps $C[0,R]\times C[0,R]$ into $C[0,R]\times
C[0,R]$. \item[(ii)] $(u,v)\in C[0,R]\times C[0,R]$ is a solution
to $(S_R)$ if and only if $T(u,v)=(u,v)$ \item[(iii)] $T$ is
completely continuous.
\end{enumerate}
See \cite{CMM93}. \item[(B)] Assume $(H_1)$. If $T(u,v)=(u,v)$,
$(u,v)\in C[0,R]\times C[0,R]$, and if there exists $\bar
r\in[0,R)$ such that $(u,v)(\bar r)\not=(0,0)$, then for any
$r\in(0,R)$ it must be that $u(r), \ v(r)>0$ and $u'(r), \
v'(r)<0$. \cite[Lemma III.1]{CFMdT} \item[(C)] Assume $(H_1)$.
Then there exists $\rho_1>0$ such that
\begin{enumerate}
\item[(i)] For any $\rho\in (0,\rho_1]$, $(0,0)$ is the only fixed
point of $T$ in $\overline{B}_\rho$. \item[(ii)] For any $\rho\in (0,\rho_1)$,  the Leray-Schauder degree
$d_{LS}(I-T,{B}(0,\rho),0)=1$. \cite[Proposition III.2]{CFMdT}
\end{enumerate}
\item[(D)] Assume $(H_1)$. Then we can choose $\eta\in
((p-1)(q-1),\delta \mu )$ and set $\theta=\eta/\delta (q-1)$. Then
$$\frac{p-1}{\delta }<\theta<\frac{\mu }{q-1},$$
and for $(u,v)\in C[0,R]\times C[0,R]$ we can define
$T_\mu(u,v)=T(u+\mu, v+\mu^\theta)$. As in \cite{CFMdT}, $T_0=T$
and $\{T_\mu\}$ is a family of compact operators satisfying the
assumptions of the homotopy theorem on any bounded interval
$[0,\bar\mu]$. If either
\begin{enumerate}
\item[(i)] there exists $\bar\mu>0$ such that $(u_\mu, v_\mu)$ is
a solution of $T_\mu(u,v)=(u,v)$ for some $\mu\ge 0$, then $\mu\le
\bar\mu$, \quad or
 \item[(ii)] there exists $M>0$ such that $(u_\mu, v_\mu)$ is
a solution of $T_\mu(u,v)=(u,v)$ with $\mu\le\bar\mu+1$, then
$||u_\mu||_\infty+||v_\mu||_\infty\le M$,
\end{enumerate}
is not satisfied, then $(S)$ has a   radially symmetric solution
$(u,v)\in C^1(\RR^N)\cap C^2(\RR^N\setminus\{0\})$ such that $u,
v$ are decreasing in $(0,\infty)$  and
$$ 0<u(r)\le Cr^{-\alpha},\quad 0<v(r)\le Cr^{-\beta},$$
where $\alpha,\ \beta$ are defined by \eqref{defsigma}.
\cite[Proposition IV.1]{CFMdT}
\end{itemize}

The proof of Proposition \ref{p1.1} is now very simple:

Since  $(S_\infty)$ does not have any nontrivial solution,  by \cite[Proposition IV.1]{CFMdT}, both $(D)(i)$ and $(D)(ii)$ must
be satisfied. Hence $d_{LS}(I-T_\mu, B(0,M+1),0)$ is well defined and has value $0$ for all $\mu\le \bar\mu+1$, and by the
homotopy theorem
\ben
d_{LS}(I-T, B(0,M+1),0)&=& d_{LS}(I-T_0, B(0,M+1),0) \\
&=&d_{LS}(I-T_{\bar\mu+1}, B(0,M+1),0)=0.
\een
Hence by $(C)(ii)$ and the excision property of the degree, there exists a
nontrivial fixed point of $T$, which by $(A)(ii)$ is a solution to $(S_R)$.

\end{proof}

\begin{remark}\label{poneeq} {\rm Let $u\ge 0$, $v\ge 0$ satisfy $(S_R)$ for $p=q=m>1$ and $\mu=\delta>0$.
Then, since from $(S_R)$, one has that  $u$ and $v$ satisfy
\ben
0&\ge& \int_B(v^\delta(x)-u^ \delta(x))(u(x)-v(x))dx\\
&=&\int_B(|\nabla u|^{m-2}\nabla u-|\nabla v|^{m-2}\nabla v)\cdot(\nabla u-\nabla v)dx\ge 0,
\een
it follows immediately that $u=v.$}
\end{remark}

\end{document}